\documentclass[12pt,leqno]{article}
\usepackage{amsmath,amssymb,amsthm,graphicx,float}

\newcommand{\bbC}{\mathbb{C}}
\newcommand{\bbD}{\mathbb{D}}
\newcommand{\bbH}{\mathbb{H}}
\newcommand{\bbN}{\mathbb{N}}

\newcommand{\go}{\omega}

\newcommand{\dcap}{\textnormal{dcap}}
\newcommand{\hcap}{\textnormal{hcap}}
\newcommand{\crad}{\textnormal{crad}}

\newcommand{\dist}{\textnormal{dist}}
\newcommand{\hyp}{\textnormal{hyp}}

\renewcommand{\Im}{\textnormal{Im}}
\newcommand{\gve}{\varepsilon}

\newcommand{\half}{\frac{1}{2}}

\newtheorem{thm}{Theorem}[section]
\newtheorem{prop}[thm]{Proposition}

\newtheorem{cor}[thm]{Corollary}

\theoremstyle{definition}
\newtheorem*{remark}{Remark}

\begin{document}

\title{Half-plane capacity and conformal radius}
\author{Steffen Rohde\footnote{Research supported by NSF Grant DMS-0800968.} ~ and Carto Wong }
\maketitle

\begin{abstract}
  In this note, we show that the half-plane capacity   of a subset of the upper half-plane is comparable to a 
  simple geometric quantity, namely the euclidean area  of  the hyperbolic neighborhood of radius one
  of this  set.  This is achieved by proving a similar estimate for the conformal radius of a subdomain of
  the unit disc, and by establishing a simple relation between these two quantities.
\end{abstract}

\section{Introduction and results}

Let $\bbH = \{ z \in \bbC \colon \Im z > 0 \}$ be the upper half-plane. A bounded subset $A \subset \bbH$
is called a {\it hull} if $\bbH\setminus A$ is a simply connected region. The {\it half-plane capacity} of a hull $A$ is the quantity
\[
    \hcap(A) := \lim_{z \to \infty} z \left[ g_A(z) - z \right],
\]
where $g_A \colon \bbH \setminus A \to \bbH$ is the unique conformal map satisfying the hydrodynamic normalization
$g(z) = z + O(\frac{1}{z})$ as $z \to \infty$. It appears frequently in connection with the Schramm-Loewner Evolution
SLE, since it serves as the conformally natural parameter in the chordal Loewner equation, see \cite{L}. In the study of
SLE, one often needs  estimates of $\hcap(A)$ in terms of geometric properties of $A$. The definition of $\hcap$ in terms
of conformal maps (or in terms of Brownian motion as in  \cite{L}) does not  immediately yield such estimates. The
purpose of this note is to provide a geometric quantity that is  comparable to  $\hcap(A)$, via a simple relation between
half-plane capacity and conformal radius.

\begin{thm} \label{T1}
  The half-plane capacity and the  (euclidean) area of the hyperbolic neighborhood of radius one are comparable,
  \[
      \hcap(A) \asymp  \left| N(A) \right|.
  \]    
\end{thm}

More precisely, there are absolute constants $C_1$, $C_2 > 0$ so that
\[
      C_1 \left| N(A) \right| \leq \hcap(A) \leq C_2 \left| N(A) \right|,
\]
where $\left| N(A) \right|$ denotes the (euclidean) area of the hyperbolic neighborhood of radius one of $A$ and
\[
    N(A) = \{ z \in \bbH \colon \dist_{\hyp}(z,A) \leq1 \}.
\]    

\begin{figure}[H]
    \begin{center}
        \includegraphics[scale=0.8]{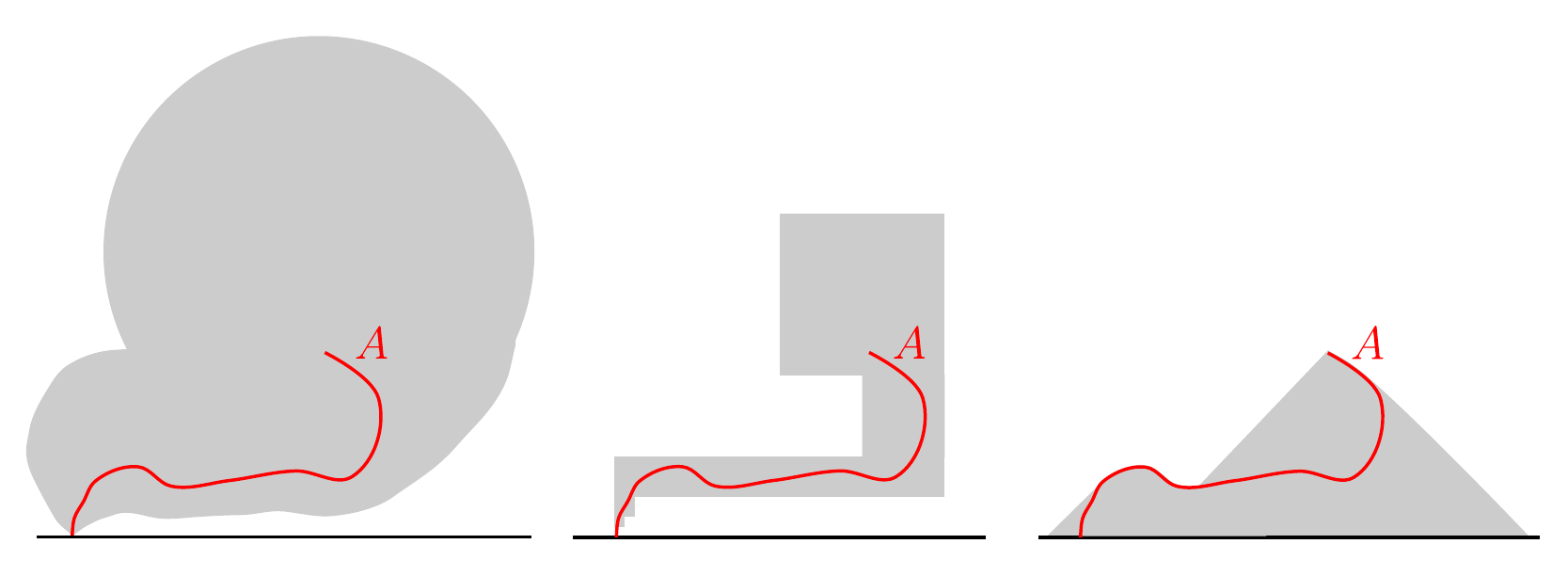}
    \end{center}
    \caption{The figure shows the hyperbolic neighborhood of $A$, the union of all Whitney squares that intersect $A$, and the area
                    under the minimal Lipschitz function of norm 1 that lies above $A$. These areas are comparable to the others with universal constants.}
                    \label{Fig:geomQuantity}
\end{figure}

Replacing the radius one by any other number only affects the constants, of course. The area of $N(A)$ is easily seen to be 
comparable to a number of other geometrically defined quantities, such as the area of all Whitney squares of $\bbH$
that intersect $A$, or the area under the minimal Lipschitz function of norm 1 that lies above $A$. See Figure~\ref{Fig:geomQuantity}.


Another important quantity is the {\it conformal radius} of a simply connected domain $D$ with respect to a point $z \in D$,
defined as
\[
    \crad(D,z) := \left| f'(0) \right|,
\]
where $f \colon \bbD \to D$ is a conformal map from the open unit disk $\bbD$ onto $D$ with $f(0)=z.$ By the Schwarz Lemma and the Koebe
One-Quarter Theorem, the conformal radius is comparable to the distance of $z$ to the boundary of $D$,
\[
    \operatorname{dist}(z,\partial D)  \leq  \crad(D,z) \leq 4 \operatorname{dist}(z,\partial D).
\]    
This geometric estimate is extremely useful in geometric function theory, but it is useless in the context of the radial Loewner equation,
where the domain $D$ is of the form $D=\bbD\setminus B$ and $B$ is small, so that $\crad(D,0)$ is close to one. In this situation,
\[
    \dcap(B) := - \log \crad(\bbD \setminus B, 0) \asymp 1 - \crad(\bbD \setminus B, 0)
\]
and we have the following theorem.

\begin{thm} \label{T2}
  If $B \subset \{ z \in \bbD \colon 1/2 < |z| <1 \}$ and if $\bbD \setminus B$ is simply connected, then 
  \[
      \dcap(B) \asymp \left| N(B) \right|.
  \]
\end{thm}

Theorems \ref{T1} and \ref{T2} are essentially equivalent. A special case of Theorem \ref{T2}, Proposition \ref{Prop1} below, 
was proven in \cite{RZ}. Until recently, Theorem~\ref{T1} did not exist in the published literature, and only existed in form
of an unpublished manuscript \cite{W}.  Recently, a probablistic proof of Theorem~\ref{T1} was given in \cite{LLN}. We give a complex
analytic proof of Theorem \ref{T2} based on \cite{RZ}. The connection between Theorems \ref{T1} and \ref{T2} is then provided by the
following simple estimates, valid for all hulls $A\subset\overline\bbH$, all bounded sets $S \subset \bbH$, and the conformal maps
$T_y(z) = \frac{z-iy}{z+iy}$ between $\bbH$ and $\bbD$ with $y > 0$:
\[
      \dcap \, T_y(A) \sim \frac2{y^2}\, \hcap(A) \qquad \text{and} \qquad
      \left| T_y(S) \right| \sim \frac{4}{y^2} \left| S \right|
\]
as $y \to \infty$, where $a \sim b$ means $\frac{a}{b} \to 1$.

\begin{remark}
  After a first version of this paper, we became aware of the papers \cite{D} and \cite{D1}. It showed \cite{D} that the half-plane capacity of a set is
  non-increasing under various notions of symmetrizations. Even though the problems and results of \cite{D} and our work are very different,
  both exploit the close relation between half-plane capacity and conformal radius. In particular, \cite[Lemma~2]{D} is equivalent to
  Corollary~\ref{Cor:hcap_1-crad} below. This allows us to define $\hcap(A)$ as a coefficient in the asymptotic expansion
  \[
      \frac{\crad(\bbH \setminus A, iy)}{\crad(\bbH, iy)} = 1 - \frac{2\, \hcap(A)}{y^2} + o(y^{-2})
  \]
  as $y \to \infty$. Under mild assumptions on a (not necessarily simply connected) domain $D \subset \widehat{\bbC}$ with an accessible boundary point
  $z_0 \in \partial D$, Dubinin and Vuorinen \cite{D1} used a similar asymptotic expansion to define the \emph{relative capacity} of a relatively closed subset
  $E \subset D$. When $(D, z_0) = (\bbH, \infty)$, relative capacity coincides with the half-plane capacity \cite[Theorem~2.6]{D1}. The behavior
  of the relative capacity under various symmetrizations and under some geometric transformations are also proved in \cite{D1}, as well as the relation
  between relative capacity and Schwarzian derivative.
\end{remark}

\section{Proofs}\label{s:proofs}

\begin{figure}[H]
  \begin{center}
    \includegraphics[scale=0.5]{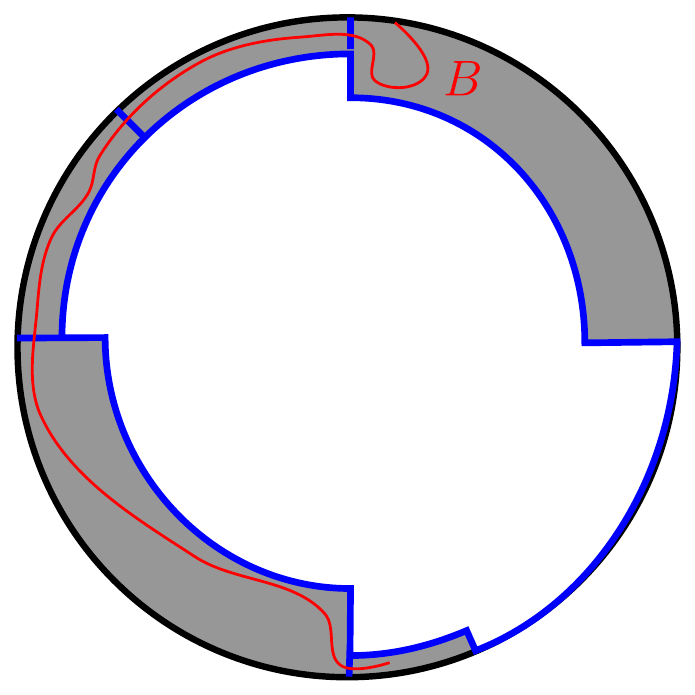}
  \end{center}
  \caption{Q(B) is the union of all dyadic squares $Q$ with $T(Q) \cap B \not= \emptyset$.} \label{Fig:Q(B)}
\end{figure}

Let  $B \subset \{ z \in \bbD \colon \half < \left| z \right| <1\}$. Denote again $N(B)$ the hyperbolic (with respect to $\bbD$) neighborhood
of radius one. Let $\widehat{N}(B)$ be the union of $N(B)$ and all of its complementary components with respect to $\bbD$ that do not
contain zero (so that $\bbD \setminus \widehat{N}(B)$ is simply connected). For a  dyadic interval $J = [\frac{k-1}{2^n}, \frac{k}{2^n}]$
($n = 1$, 2, \dots and $k = 1$, $2$, \dots,~$2^n$), consider the dyadic ``square''
\[
    Q_J = \left\{ z  \in \bbD \colon \frac{z}{\left|z\right|} \in e^{2 \pi i J}, 1 - \left| z \right| \leq \frac{1}{2^n}  \right\}
\]    
and its top half $T(Q_J) = \{ z \in Q_J \colon 1 - \left|z\right| > \frac{1}{2^{n+1}} \}$. Denote $Q(B)$ the union of all dyadic squares $Q$
with $T(Q) \cap B \not= \emptyset$, see Figure~\ref{Fig:Q(B)}. Clearly, $\left|Q(B)\right|$, $\left| N(B) \right|$  and $\left|\widehat{N}(B)\right|$
are comparable with universal constants. The proof of Proposition 2 in \cite{RZ} showed the following.

\begin{prop}\label{Prop1}
  If $B \subset \{ z \in \bbD \colon \half < \left| z \right| < 1 \}$ and $\bbD\setminus B$ is a simply connected region, then
  \[ 
      C_1 \left| B \right| \leq \dcap(B) \leq \dcap \, Q(B) \leq C_2 \left| Q(B) \right|,
  \]
  with absolute constants $C_1$, $C_2 > 0$.
\end{prop}

\begin{figure}[H]
    \begin{center}
        \includegraphics[scale=0.9]{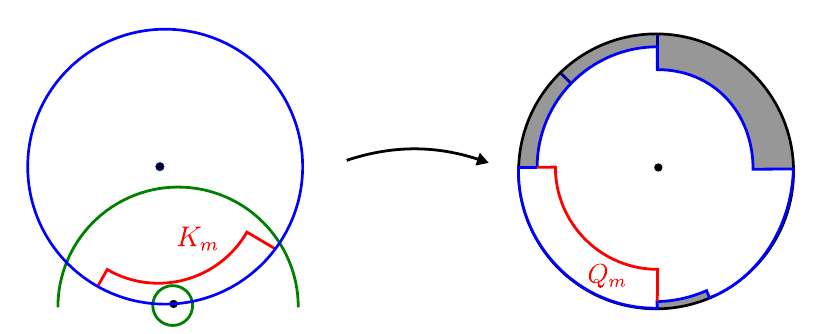}
    \end{center}
    \caption{A proof of (\ref{I:induction}) uses the fact that $\{K_m\}$ is a family of uniform quasi-circles.} \label{Fig:induction}    
\end{figure}

 \begin{proof}[Sketch of proof]
    The first inequality in Proposition~\ref{Prop1} follows from Parseval's formula. The second inequality is trivial (Schwarz's Lemma). The
    final inequality is proved  by an inductive procedure as follows. Write
    \[
        Q(B) = \bigcup_{j=1}^{\infty} Q_j,
    \]
    where $\{ Q_j \}$ is a disjoint (modulo boundary) family of dyadic squares, arranged so that $\left| Q_1 \right| \geq \left| Q_2 \right| \geq \cdots$.
    It suffices to show that
    \begin{equation} \label{I:induction}
        \dcap \left( \bigcup_{j=m}^{\infty} Q_j \right) - \dcap \left( \bigcup_{j=m+1}^{\infty} Q_j \right) \asymp \left| Q_m \right|
    \end{equation}
    for every $m$. Let $f_m$ be the conformal map from $\bbD$ onto $\bbD \setminus \bigcup_{j=m+1}^{\infty} Q_j$, normalized so that $f_m(0) = 0$.
    Let $K_m = f_m^{-1}(Q_m)$. Inequality (\ref{I:induction}) is equivalent to $\dcap(K_m) \asymp \left| Q_m \right|$. In \cite{RZ}, this last inequality
    was proved by noting that $\{ K_m \}$ is a family of uniform quasi-circles. This allows us to construct two concentric circular hulls of comparable sizes,
    so that one contains $K_m$ and the other one is contained in $K_m$. See Figure~\ref{Fig:induction}.
\end{proof}

The proof of Theorem \ref{T2} follows from Proposition \ref{Prop1},  the comparability of $\left|Q(B)
\right|$  and $\left| N(B) \right|$, and the fact that ``fattening'' the set $B$ does not change $\dcap(B)$ by
more than a factor:

\begin{prop} \label{Prop:fattening}
  If $B \subset \{ z \in \bbD \colon \frac{1}{2} < \left| z \right| < 1 \}$ and $\bbD \setminus B$ is
  a simply connected region, then
  \[
      \dcap \, \widehat{N}(B) \leq C \, \dcap(B),
  \]
  where $C > 0$ is an absolute constant.  
\end{prop}

\begin{figure} \label{Fig:harmEst}
    \begin{center}
        \includegraphics[scale=0.8]{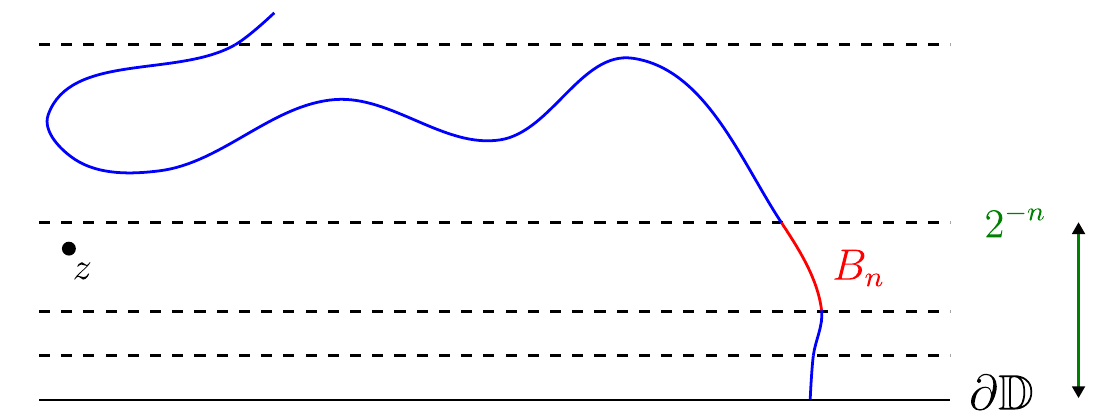}
    \end{center}    
    \caption{In the proof of Proposition~\ref{Prop:fattening}, it is not generally true that $\widetilde{\go}_n(z) \lesssim \go_n(z)$.
    The figure illustrates a situation where $\widetilde{\go}(z) = 1$ but $\go_n(z)$ is small.}
\end{figure}

\begin{proof}
  Let $\gve > 0$ be a small number that will be chosen later. Define $\widehat{N}_{\gve}(B)$ the same way as $\widehat{N}(B)$ using
  hyperbolic radius $\gve$ instead of one. 
  We only need to show
  \begin{equation} \label{I:1-crad}
      \dcap \, \widehat{N}_{\gve}(B) \lesssim \dcap(B),
  \end{equation}
  where $a \lesssim b$ means $a \leq C b$ for some constant $C > 0$. Once we have proved (\ref{I:1-crad}), iterating the inequality
  $\frac{1}{\gve}$ times gives $\dcap \, \widehat{N}(B) \lesssim \dcap(B)$. The reverse inequality is
  trivial by the Schwarz's Lemma.
   
  Let $\widehat{B} = \widehat{N}_{\gve}(B)$.  Denote $f \colon \bbD \to \bbD \setminus B$ the conformal map under the normalization
  $f(0) = 0$ and $f'(0) > 0$. Since $z \mapsto \log \left| \frac{f(z)}{z} \right|$ is harmonic on $\bbD$, the mean value property gives
  \begin{equation} \label{E:log|f'(0)|}
      \log \left| f'(0) \right| = \frac{1}{2 \pi} \int_0^{2 \pi} \log \left| f(e^{it}) \right| \, dt.
  \end{equation}
  (If $f$ does not continuously extend to $\partial \bbD$, either approximate $\bbD \setminus B$ by
  smooth regions, or interpret $f(e^{it})$ as an angular limit.) We decompose $\bbD$ into an union of dyadic layers.
  For each $n = 1$, 2, \dots, let
  \[
      D_n = \left\{ z \in \bbD \colon \frac{1}{2^{n+1}} \leq 1 - \left|z\right| < \frac{1}{2^n} \right\} \qquad
      \text{and} \qquad B_n = B \cap D_n.
  \]
  From (\ref{E:log|f'(0)|}) and the elementary fact that $1 - u \asymp - \log(u)$  for $\frac{1}{2} \leq u < 1$, we have
  \begin{equation} \label{E:crad_1}
      \dcap(B) \asymp \sum_{n=1}^{\infty} \frac{\go_n(0)}{2^n},
  \end{equation}
  where $\go_n(z) = \go(z, B_n, \bbD \setminus B)$ denotes the harmonic measure of $B_n$ with respect to the
  region $\bbD \setminus B$ and the point $z$. Define $\widehat{\go}_n(z)$ similarly using $\widehat{B}$ instead of $B$. We have
  \begin{equation} \label{E:crad_2}
      \dcap(\widehat{B}) \asymp \sum_{n=1}^{\infty} \frac{\widehat{\go}_n(0)}{2^n}.
  \end{equation}
  Recall that our $\widehat{B} = \widehat{N}_{\gve}(B)$ depends on $\gve$. Choose $\gve > 0$ so that for each $n \in \bbN$,
  every hyperbolic ball (in $\bbD$) centered in $D_n$ of radius $2 \gve$ is contained in $D_{n-1} \cup D_n \cup D_{n+1}$.
  This is possible because the hyperbolic distance $\dist_{\text{hyp}}(D_n, D_{n+2}) \asymp 1$. We claim that
  \begin{equation} \label{I:w_n(z)}
      \widehat{\go}_n(z) \lesssim \go_{n-1}(z) + \go_n(z) + \go_{n+1}(z)
  \end{equation}
  for all $n \in \bbN$ and all $z \in \bbD \setminus \widehat{B}$. (The inequality $\widehat{\go}_n(z) \lesssim
  \go_n(z)$ is not necessarily true, see Figure~\ref{Fig:harmEst}.) If we can show (\ref{I:w_n(z)}), then (\ref{E:crad_1}), (\ref{E:crad_2}),
  (\ref{I:w_n(z)}) imply (\ref{I:1-crad}) and complete the proof.
  
  Since both sides of (\ref{I:w_n(z)}) are harmonic functions on $\bbD \setminus \widehat{B}$, by the maximum principle it suffices
  to show (\ref{I:w_n(z)}) for $z$ in the boundary of this region. If $z \in \partial (\bbD \setminus \widehat{B})$ and $z \notin D_n$,
  then $\widehat{\go}_n(z) = 0$ and (\ref{I:w_n(z)}) is trivial in this case. Suppose $z \in \partial (\bbD \setminus \widehat{B})$ and
  $z \in D_n$ from now on. Since $\widehat{\go}_n(z) \leq 1$, all we need to show is that the harmonic measure
  \begin{equation} \label{Q:harm}
      \go(z, B_{n-1} \cup B_n \cup B_{n+1}, \bbD \setminus B)
  \end{equation}
  is bounded away from zero. Our choice of $\gve$ guarantees that the hyperbolic ball $B_{\text{hyp}}(z,2\gve)$ centered at $z$
  with radius $2 \gve$ is contained in $D_{n-1} \cup D_n \cup D_{n+1}$, so that the set $E := B_{n-1} \cup B_n \cup B_{n+1}$ connects
  the hyperbolic circles $\partial B_{\text{hyp}}(z, \gve)$ and $\partial B_{\text{hyp}}(z, 2\gve)$. Now, the maximum principle and
  the Beurling projection theorem provide a lower bound of (\ref{Q:harm}):
  \[
      \go(z, E, \bbD \setminus B) \geq \go(z, E, B_{\text{hyp}}(z,2\gve)) \geq \go(0, [r,1], \bbD)
  \]
  for some absolute constant $r \in (0,1)$.
\end{proof}

\begin{proof}[Proof of Theorem~\ref{T2}]
  By Proposition~\ref{Prop:fattening} and Proposition~\ref{Prop1},
  \[
      \dcap(B) \gtrsim \dcap \, \widehat{N}(B) \gtrsim \left| \widehat{N}(B) \right|
      \geq \left| N(B) \right|.
  \]
  On the other hand,
  \[
      \begin{aligned}
          \dcap(B) &\leq \dcap \, Q(\widehat{N}(B)) \\
          &\lesssim \left| Q(\widehat{N}(B)) \right| \\
          &\lesssim \left| N(B) \right|
      \end{aligned}
  \]
  by the Schwarz's Lemma, Proposition~\ref{Prop1} and the simple facts that $Q(\widehat{N}(B)) = Q(N(B))$ and $\left| Q(N(B)) \right|
  \lesssim \left| N(B) \right|$.
\end{proof}

We now prove Theorem~\ref{T1} using Theorem~\ref{T2}. In fact, the two theorems are equivalent. Proposition~\ref{Prop:hcap_crad}
below gives a connection between half-plane capacity and conformal radius. Roughly speaking, the half-plane capacity of a hull is
(up to a constant) asymptotically equal to the change of conformal radius, with respect to a point close to infinity, when the hull is
removed. By scaling, we may take the reference point to be $i$ and consider hulls that are small.

\begin{prop} \label{Prop:hcap_crad}
  There are absolute constants $C > 0$ and $\gve_0 \in (0,1)$ such that for any hull $A \subset \bbH$ with $\sup_{z \in A} \left| z \right|
  \leq \gve \leq \gve_0$,
  \begin{equation} \label{I:hcap_crad}
      \left| \frac{2 - \crad(\bbH \setminus A, i)}{\hcap(A)} - 4 \right| \leq C \gve.
  \end{equation}
\end{prop}

For our purpose, we will only need the weaker statement that
\[
    \lim_{y \to \infty} \frac{2 - \crad(\bbH \setminus (y^{-1} A), i)}{\hcap(y^{-1} A)} = 4
\]
for every hull $A \subset \bbH$.

\begin{proof}
  Note that $\crad(\bbH,z_0) = 2 \, \Im(z_0)$. We have
  \[
      \crad(\bbH \setminus A, i) = \frac{2 \, \Im g(i)}{\left| g'(i) \right|},
  \]
  where $g = g_A \colon \bbH \setminus A \to \bbH$ is the hydrodynamically normalized conformal map. Write $h = \hcap(A)$ for the sake
  of notation. We claim that
  \begin{align}
          \left| g(i) - i + ih \right| &\leq C_1 h \gve \label{I:g(i)} \\
           \left| \frac{1}{g'(i)} - 1 + h \right| &\leq C_2 h \gve \label{I:g'(i)}
  \end{align}
  for some absolute constants $C_1$, $C_2 > 0$. These two inequalities imply
  \[
      \left| \frac{\Im g(i)}{g'(i)} - (1-h)^2 \right| \leq C_3 h \gve,
  \]
  and (\ref{I:hcap_crad}) follows.

  To prove (\ref{I:g(i)}), we use the Poisson integral representation
  \[
      g(z) - z = \frac{1}{\pi} \int_{-2\gve}^{2\gve}  \frac{\Im[g^{-1}(x)]}{g(z)-x} \, dx
      \qquad (z \in \bbH \setminus A)
  \]
  which implies that $h = \lim_{z \to \infty} z [g(z) - z] = \frac{1}{\pi} \int_{-2\gve}^{2\gve}  \Im[g^{-1}(x)] \, dx$ and therefore
  \begin{equation} \label{E:g(z)-z-h/z}
      g(z) - z - \frac{h}{z} = \frac{1}{\pi} \int_{-2\gve}^{2\gve} \Im[g^{-1}(x)] \left( \frac{1}{g(z)-x} - \frac{1}{z} \right) \, dx
      \qquad (z \in \bbH \setminus A).
  \end{equation}
  On the other hand, we have $\left| g(z) - z \right| \leq 3 \gve$. (See \cite[Corollary~3.44]{L}, which showed this estimate
  using the maximum principle and the observation that it holds for $z \in \partial (\bbH \setminus A)$.) Using this simple estimate,
  we can prove that for all $x \in [-2\gve, 2\gve]$,
  \begin{equation} \label{I:integralKernel}
      \left| \frac{1}{g(z)-x} - \frac{1}{z} \right| \leq \frac{5 \gve}{\left| z \right| (\left| z \right| - 5 \gve)}.
  \end{equation}
  Inequality (\ref{I:g(i)}) follows from (\ref{E:g(z)-z-h/z}) and (\ref{I:integralKernel}) with $z = i$. Finally, (\ref{I:g'(i)}) can be proved similarly
  using
  \[
      \frac{z_2 - z_1}{g(z_2) - g(z_1)} - 1 - \frac{h}{z_1 z_2} = \frac{1}{\pi} \int_{-2\gve}^{2\gve} \Im[g^{-1}(x)]
      \left( \frac{1}{(g(z_1)-x)(g(z_2)-x)} - \frac{1}{z_1 z_2} \right) \, dx
  \]
  instead of (\ref{E:g(z)-z-h/z}).
\end{proof}

\begin{cor} \label{Cor:hcap_1-crad}
  If $A \subset \bbH$ is a hull and $B_y = T_y(A) \subset \bbD$, where $T_y(z) = \frac{z - iy}{z + iy}$ and $y > 0$, then
  \[
      \lim_{y \to \infty} \frac{y^2 \dcap(B_y)}{\hcap(A)} = 2.
  \]
  In other words, $\dcap(B_y) \sim \frac{2}{y^2} \hcap(A)$ as $y \to \infty$.
\end{cor}

\begin{proof}[Proof of Theorem~\ref{T1}]
  Fix a hull $A \subset \bbH$ and let $T_y(z) = \frac{z - iy}{z + iy}$ as before. When $y > 0$ is large, the hulls $\frac{1}{y} A$ and
   $B_y = T_y(A) = T_1(y^{-1} A)$ are both small. By Corollary~\ref{Cor:hcap_1-crad} and Theorem~\ref{T2},
  \[
      \frac{1}{y^2} \hcap(A) \asymp \dcap(B_y) \asymp \left| N(B_y) \right| \asymp \left| N(y^{-1}A) \right|
      \asymp \frac{1}{y^2} \left| N(A) \right|.
  \]
\end{proof}


\end{document}